\documentstyle[amsmath,amsfonts,amssymb,amsthm,12pt]{article}
\font\sixbb=msbm6
\font\eightbb=msbm8
\font\twelvebb=msbm10 scaled 1095
\newfam\bbfam
\textfont\bbfam=\twelvebb \scriptfont\bbfam=\eightbb
                           \scriptscriptfont\bbfam=\sixbb
\def\bb{\fam\bbfam\twelvebb}
\newcommand{\Rea}{{\bb R}}

\newcommand{\KK}{{\bb K}}

\newtheorem{theorem}{\bf Theorem}[section]
\newtheorem{claim}[theorem]{\bf Claim}

\newtheorem{proposition}[theorem]{\bf Proposition}

\newcommand{\enp}{\begin{flushright} $\Box$ \end{flushright}}
\newcommand{\beq}[0]{\begin{equation}}
\newcommand{\enq}[0]{\end{equation}}
\newcommand{\th}{\tilde{{\rm H}}}

\newcommand{\cf}{{\cal F}}
\newcommand{\reg}{{\rm reg}}
\newcommand{\pd}{{\rm pd}}
\newcommand{\ck}{{\cal K}}

\newcommand{\cl}{{\cal L}}

\newcommand{\St}{{\rm St}}
\newcommand{\lk}{{\rm lk}}
\newcommand{\sd}{{\rm sd}}

\newcommand{\dpo}{\overset{.}{D}}
\newcommand{\hell}{{\rm h}}

\newcommand{\hv}{{\rm H}}
\newcommand{\tho}{\tilde{\rm H}}
\newcommand{\dh}{{\rm h}}
\newcommand{\dth}{\tilde{\rm h}}
\newcommand{\LK}{{\rm L}_{\KK}}

\title{Intersections of Leray Complexes and \\ Regularity of Monomial Ideals}
\author{Gil Kalai\thanks{Institute of Mathematics, Hebrew
University, Jerusalem 91904, Israel,  and Departments of Computer
Science and Mathematics, Yale University. e-mail:
kalai@math.huji.ac.il .  Research supported by ISF, BSF and NSF
grants.} \and Roy Meshulam\thanks{Department of Mathematics,
Technion, Haifa 32000, Israel,  and Institute for Advanced Study,
Princeton 08540. e-mail: meshulam@math.technion.ac.il~. Research
supported by the Israel Science Foundation and by a State of New
Jersey grant. }}
\date{}
\begin{document}
\insert\footins{\footnotesize\rule{0pt}{\footnotesep}
\\ {\it Math Subject Classification.}  13F55, 55U10.
\\ }
\maketitle
\begin{abstract}
For a simplicial complex $X$ and a field $\KK$, let
$\dth_i(X)=\dim \tho_i(X;\KK).$ \\  It is shown that if $X,Y$ are
complexes on the same vertex set,  then for $k \geq 0$
$$
\dth_{k-1}(X\cap Y) \leq \sum_{\sigma \in Y} \sum_{i+j=k}
\dth_{i-1}(X[\sigma])\cdot \dth_{j-1}(\lk(Y,\sigma))~~.
$$
$~~$ A simplicial complex $X$ is {\it $d$-Leray} over $\KK$, if
$\tho_i(Y;\KK)=0$ for all induced subcomplexes $Y \subset X$ and
$i \geq d$. Let $L_{\KK}(X)$ denote the minimal $d$ such that $X$
is $d$-Leray over $\KK$.  The above theorem implies that if $X,Y$
are simplicial complexes on the same vertex set then
$$L_{\KK}(X \cap Y) \leq L_{\KK}(X) + L_{\KK}(Y) ~~~.$$
$~$ Reformulating this inequality in commutative algebra terms, we
obtain the following result conjectured by Terai: If $I,J$ are
square-free monomial ideals in $S=\KK[x_1,\ldots,x_n]$, then
$$\reg(I+J) \leq \reg(I)+\reg(J)-1~~$$ where $\reg(I)$ denotes the
Castelnuovo-Mumford regularity of $I$.

\end{abstract}

\section{Introduction}

$~~~$ Let $X$ be a simplicial complex on the vertex set $V$. The
{\it induced} subcomplex on a subset of vertices $S \subset V$ is
$X[S]=\{ \sigma \in X: \sigma \subset S \}$. Let $\{\}$ be the
{\it void complex} and let $\{\emptyset\}$ be the {\it empty
complex}. Any non-void complex contains $\emptyset$ as a unique
$(-1)$-dimensional face. The {\it star} of a subset $A \subset V$
is $\St(X,A)=\{\tau \in X:\tau \cup A \in X ~\}~$. The {\it link}
of $A \subset V$ is $\lk(X,A)=\{\tau \in \St(X,A):\tau \cap A
=\emptyset~\}~.$
 If $A \not\in X$ then
$\St(X,A)=\lk(X,A)=\{\}$. All homology groups considered below are
with coefficients in a fixed field $\KK$ and we denote
$\dth_i(X)=\dim_{\KK} \tho_i(X)$. Note that $\dth_{-1}(\{\})=0
\neq 1 = \dth_{-1}(\{\emptyset\})$. Our main result is the
following
\begin{theorem}
\label{dimt} Let $X,Y$ be finite simplicial complexes on the
same vertex set. Then for $k \geq 0$
\begin{equation}
\label{eq2} \dth_{k-1}(X\cap Y) \leq \sum_{\sigma \in Y}
\sum_{i+j=k} \dth_{i-1}(X[\sigma])\cdot
\dth_{j-1}(\lk(Y,\sigma))~~.
\end{equation}
\end{theorem}
We next discuss some applications of Theorem \ref{dimt}. A
simplicial complex $X$ is {\it $d$-Leray} over $\KK$ if
$\tho_i(Y)=0$ for all induced subcomplexes $Y \subset X$ and $i
\geq d$. Let $\LK (X)$ denote the minimal $d$ such that $X$ is
$d$-Leray over $\KK$. Note that $\LK(X)=0$ iff $X$ is a simplex.
$\LK(X) \leq 1$ iff $X$ is the clique complex of a chordal graph
(see e.g. \cite{W75}).

The class $\cl_{\KK}^d$ of $d$-Leray complexes over $\KK$ arises
naturally in the context of  Helly type theorems \cite{E93}. The
{\it Helly number} $~\hell (\cf)$ of a finite family of sets $\cf$
is the minimal positive integer $h$ such that if $\ck \subset \cf$
satisfies $\bigcap_{K \in \ck'}K \neq \emptyset$ for all $\ck'
\subset \ck$ of cardinality $\leq h$, then $\bigcap_{K \in \ck} K
\neq \emptyset$.  The {\it nerve} $N(\ck)$ of a family of sets
$\ck$, is the simplicial complex whose vertex set is $\ck$ and
whose simplices are all $\ck' \subset \ck$ such that $ \bigcap_{K
\in \ck'} K \neq \emptyset$. It is easy to see that for any field
$\KK$
$$\hell (\cf) \leq 1+ \LK(N(\cf)).$$ For example, if $\cf$ is a
finite family of convex sets in $\Rea^d$, then  by the Nerve Lemma
(see e.g. \cite{B95}) $N(\cf)$ is $d$-Leray over $\KK$, hence
follows Helly's Theorem: $\hell(\cf) \leq d+1$. This argument
actually proves the Topological Helly Theorem: If $\cf$ is a
finite family of closed sets in $\Rea^d$ such that the
intersection of any subfamily of $\cf$ is either empty or
contractible, then $\hell(\cf) \leq d+1$.

Nerves of families of convex sets however satisfy a stronger
combinatorial property called {\it $d$-collapsibility} \cite
{W75}, that leads to some of the deeper extensions of Helly's
Theorem. It is of considerable interest to understand which
combinatorial properties of nerves of families of convex sets in
$\Rea^d$ extend to arbitrary $d$-Leray complexes. For some recent
work in this direction see \cite{akmm,KM}. One consequence of
Theorem \ref{dimt} is the following
\begin{theorem}
\label{leray} Let $X_1,\ldots,X_r$ be simplicial complexes on the
same finite vertex set. Then
\begin{equation}
\label{lerayi} \LK \bigl( \bigcap_{i=1}^r X_i \bigr) \leq
\sum_{i=1}^r \LK(X_i)
\end{equation}
\begin{equation}
\label{lerayu} \LK \bigl( \bigcup_{i=1}^r X_i \bigr) \leq
\sum_{i=1}^r \LK (X_i)+r-1 ~~~.
\end{equation}
\end{theorem}
\noindent {\bf Example:} Let $V_1,\ldots,V_r$ be disjoint sets of
cardinalities  $|V_i|=a_i$, and let $V=\bigcup_{i=1}^r V_i$. Let
$\Delta(A)$ denote the simplex on vertex set $A$, with boundary
$\partial \Delta(A)\simeq S^{|A|-2}$ . Consider the complexes
$$X_i= \Delta(V_1)* \cdots *\Delta(V_{i-1})*\partial
\Delta(V_i)*\Delta(V_{i+1})* \cdots *\Delta(V_r)~~.$$ Then
$$\bigcap_{i=1}^r X_i=\partial \Delta(V_1) * \cdots *
\partial\Delta(V_r) \simeq S^{\sum_{i=1}^ra_i-r-1}~~~$$
and
$$\bigcup_{i=1}^r X_i= \partial\Delta(V_1 \cup \ldots \cup
V_r) \simeq S^{\sum_{i=1}^r a_i -2}~~.$$ The only non-contractible
induced subcomplex of $X_i$ is $\partial\Delta(V_i)$, therefore
$\LK(X_i)=a_i-1~$. Similar considerations show that
$~\LK(\cup_{i=1}^r X_i)=\sum_{i=1}^r a_i-1~$ and
$~\LK(\cap_{i=1}^r X_i)=\sum_{i=1}^r a_i-r~$, so equality is
attained in both (\ref{lerayi}) and
(\ref{lerayu}). \\

Theorem \ref{leray} was first conjectured in a different but
equivalent form by Terai \cite{T99}, in the context of monomial
ideals . Let $S=\KK[x_1,\ldots,x_n]$ and let $M$ be a graded
$S$-module. Let $\beta_{ij}(M)=\dim_{\KK} {\rm Tor}_i^S(\KK,M)_j$
be the graded Betti numbers of $M$. The {\it regularity} of $M$ is
the minimal $\rho=\reg(M)$ such that $\beta_{ij}(M)$ vanish for
$j>i+\rho~~$ (see e.g. \cite{Eisenbud}).
\\
For a simplicial complex $X$ on $[n]=\{1,\ldots,n\}$ let $I_X$
denote the ideal of $S$ generated by $\{\prod_{i \in A} x_i: A
\not\in X\}$.
The following fundamental result of Hochster
relates the Betti numbers of $I_X$ to the topology of the induced
subcomplexes $X$.
\begin{theorem}[Hochster \cite{H75}]
\label{hochster}
\begin{equation}
\label{hoch1} \beta_{ij}(I_X)=\sum_{|W|=j} \dim_{\KK}
\tho_{j-i-2}(X[W])~~.
\end{equation}
\end{theorem}
\noindent Hochster's formula (\ref{hoch1}) implies that $\reg(I_X)=\LK(X)+1$. The case $r=2$ of Theorem
\ref{leray} is therefore equivalent to the following result
conjectured by Terai \cite{T99}.
\begin{theorem}
\label{mono}
Let $X$ and $Y$ be simplicial complexes on the same vertex set.
Then
$$\reg(I_X+I_Y) =\reg(I_{X \cap Y}) \leq \reg(I_X)+\reg(I_Y)-1~~$$
$$\reg(I_X \cap I_Y)=\reg(I_{X \cup Y}) \leq \reg(I_X)+\reg(I_Y)~~.$$
{\enp}
\end{theorem}
Theorem \ref{mono} can also be formulated in terms of projective
dimension. Let $X^*=\{\tau \subset [n]: [n]-\tau \not\in X\}$
denote the Alexander dual of $X$. Terai \cite{T99a} showed that
\begin{equation}
\label{terai} \pd(S/I_X)=\reg(I_{X^*})~~.
\end{equation}
Using (\ref{terai}) it is straightforward to check
that Theorem \ref{mono} is equivalent to
\begin{theorem}
\label{proj}
$$\pd(I_X \cap I_Y) \leq \pd(I_X)+\pd(I_Y)~~$$
$$\pd(I_X + I_Y) \leq \pd(I_X)+\pd(I_Y)+1~~.$$
{\enp}
\end{theorem}
In Section \ref{s:spec} we give a spectral sequence for the
relative homology group $\hv_*(Y,X \cap Y)$, which directly
implies Theorem \ref{dimt}. The proof of Theorem \ref{leray} is
given in Section \ref{s:leray}.

\section{A Spectral Sequence for $\hv_*(Y,X \cap Y)$}
\label{s:spec} $~~~$ Let $K$ be a simplicial complex. The
subdivision $\sd(K)$ is the order complex of the set of the
non-empty simplices of $K$ ordered by inclusion. For $\sigma \in
K$ let $D_K(\sigma)$ denote the order complex of the interval
$[\sigma,\cdot]=\{\tau \in K ~:~ \tau \supset \sigma \}$.
$D_K(\sigma)$ is called the {\it dual cell} of $\sigma$. Let
$\dpo_K(\sigma)$ denote the order complex of the interval
$(\sigma,\cdot]=\{\tau \in K ~:~ \tau \supsetneqq \sigma\}$. Note
that $\dpo_K(\sigma)$ is isomorphic to $\sd(\lk(K,\sigma))$ via
the simplicial map $\tau \rightarrow \tau-\sigma$. Since
$D_K(\sigma)$ is contractible, it follows that
$\hv_i(D_K(\sigma),\dpo_K(\sigma)) \cong \th_{i-1}(\lk(K,\sigma))$
for all $i \geq 0$. Write $K(p)$ for the family of $p$-dimensional
simplices in $K$.  The proof of Theorem \ref{dimt} depends on the
following
\begin{proposition}
\label{spec} Let $X$ and $Y$ be two complexes on the same vertex
set $V$, such that $\dim Y=n$. Then there exists a homology
spectral sequence $\{E_{p,q}^r\}$ converging to $\hv_*(Y,X \cap
Y)$ such that
$$E_{p,q}^1=\bigoplus_{\sigma \in Y(n-p)}
\bigoplus_{\substack{i,j \geq 0 \\
i+j=p+q}} \th_{i-1}(X[\sigma]) \otimes
\th_{j-1}(\lk(Y,\sigma))~~$$ for $~0 \leq p \leq n~$,$~0 \leq q~$,
and $E_{p,q}^1=0$ otherwise.
\end{proposition}
\noindent {\bf Proof:} In the sequel we identify abstract
complexes with their geometric realizations. Let $\Delta$ denote
the simplex on $V$. For $0 \leq p \leq n$ let
$$K_p=\bigcup_
{\substack{\sigma \in Y \\
\dim \sigma \geq n-p}} \Delta[\sigma] \times
D_Y(\sigma) \subset Y \times \sd(Y)~~$$ and
$$L_p=\bigcup_
{\substack{\sigma \in Y \\
\dim \sigma \geq n-p}} X[\sigma] \times D_Y(\sigma) \subset (X
\cap Y) \times \sd(Y)~~.$$ Write $K=K_n~,~L=L_n$. Let $$\pi:K
\rightarrow \bigcup_{\sigma \in Y} \Delta[\sigma] = Y~~$$ denote
the projection on the first coordinate. For a point $z \in Y$, let
$\tau={\rm supp}(z)$ denote the minimal simplex in $Y$ containing
$z$. The fiber $\pi^{-1}(z)=\{z\} \times D_Y(\tau)$ is a cone,
hence $\pi$ is a homotopy equivalence. Similarly, the restriction
$$\pi_{|L}: L \rightarrow \bigcup_{\sigma \in Y} X[\sigma] = X \cap
Y~~$$ is a homotopy equivalence. Let $F_p=C_*(K_p,L_p)$ be the
group of cellular chains of the pair $(K_p,L_p)$. The filtration
$0 \subset F_0  \subset \cdots \subset F_n=C_*(K,L)$ gives rise to
a homology spectral sequence $\{E^r\}$ converging to $\hv_*(K,L)
\cong \hv_*(Y,X \cap Y)$. We compute $E^1$ by excision and the
K\"{u}nneth formula:
$$E_{p,q}^1=\hv_{p+q}(F_p/F_{p-1}) \cong \hv_{p+q} (K_p,L_p \cup K_{p-1}) \cong
$$ $$\hv_{p+q} \bigl( \bigcup_{\sigma \in
Y(n-p)} \Delta[\sigma] \times D_Y(\sigma), \bigcup_{\sigma \in Y(n-p)}
X[\sigma] \times D_Y(\sigma)~ \cup ~\Delta[\sigma] \times \dpo_Y(\sigma) ~\bigr) ~
 \cong $$ $$ \bigoplus_{\sigma \in Y(n-p)} \hv_{p+q}\bigl(\Delta[\sigma] \times D_Y(\sigma),
X[\sigma] \times D_Y(\sigma)~ \cup ~\Delta[\sigma] \times \dpo_Y(\sigma) \bigr) \cong$$
$$  \bigoplus_{\sigma \in Y(n-p)} \bigoplus_{i+j=p+q} \hv_i(\Delta[\sigma],X[\sigma]) \otimes
\hv_j (D_Y(\sigma),\dpo_Y(\sigma)) \cong$$
$$
\bigoplus_{\sigma \in Y(n-p)} \bigoplus_{i+j=p+q}
\th_{i-1}(X[\sigma]) \otimes \th_{j-1}( \lk(Y,\sigma))~~.$$ {\enp}
\noindent {\bf Remark:} The derivation of the above spectral sequence
may be viewed as a simple application of the method of simplicial resolutions.
See Vassiliev's papers \cite{V99,V01} for a description of this technique,
and for far reaching applications to plane arrangements and to
spaces of Hermitian operators.
\ \\ \\
\noindent {\bf Proof of Theorem \ref{dimt}:} By Proposition
\ref{spec}
$$\dth_{k-1}(X \cap Y) \leq \dth_{k-1}(Y)+\dh_k(Y,X \cap Y) \leq
$$ $$\dth_{k-1}(Y)+\sum_{p+q=k} \dim E_{p,q}^1 =$$ $$
\dth_{k-1}(Y)+\sum_{\substack{\emptyset \neq \sigma \in Y \\
\dim \sigma \geq n-k}} \sum_{i+j=k}   \dth_{i-1}(X[\sigma])\cdot
\dth_{j-1}(\lk(Y,\sigma)) \leq
$$
$$ \sum_{\sigma \in Y} \sum_{i+j=k} \dth_{i-1}(X[\sigma])\cdot
\dth_{j-1}(\lk(Y,\sigma))~~.$$ {\enp}

\section{Intersection of Leray Complexes}
\label{s:leray} $~~~$ We first recall a well-known
characterization of $d$-Leray complexes. For completeness we
include a proof.
\begin{proposition}
\label{folk} For a simplicial complex $X$, the following
conditions are equivalent:
\begin{itemize}
\item[(i)] $X$ is $d$-Leray  over $\KK$.
\item[(ii)] $\th_i(\lk(X,\sigma))=0$ for every $\sigma \in X$ and $i \geq d$.
\end{itemize}
\end{proposition}
\noindent It will be convenient to prove a slightly more general
result. Let $k,m \geq 0$. We say that a simplicial complex $X$ on
$V$ satisfies {\it condition $P(k,m)$} if $\th_i \bigl(
\lk(X[A],B) \bigr) =0$ for all $B \subset A \subset V$ such that
$|A| \geq |V|-k~,~|B| \leq m$.
\begin{claim}
\label{st-lk} If $k \geq 0$ and $ m \geq 1$ then conditions
$P(k,m)$ and $P(k+1,m-1)$ are equivalent.
\end{claim}
\noindent {\bf Proof:} Suppose $B \subset A \subset V$ and $B_1
\subset A_1 \subset V$ satisfy $B=B_1 \cup \{v\}~,~A=A_1 \cup
\{v\}$ for some $v \not\in A_1$, and let
$$Z_1=\lk(X[A_1],B_1)~~~,~~~Z_2=\St(\lk(X[A],B_1),v)~~.$$
Then $$Z_1 \cup Z_2= \lk(X[A],B_1)~~~,~~~Z_1 \cap Z_2=
\lk(X[A],B)~~$$ and by  Mayer-Vietoris there is an exact sequence
$$
\ldots  \rightarrow \th_{i+1} \bigl( \lk(X[A],B_1) \bigr)
\rightarrow \th_i \bigl( \lk(X[A],B) \bigr) \rightarrow
$$
\begin{equation}
\label{mayer} \th_i \bigl( \lk(X[A_1],B_1) \bigr) \rightarrow
\th_i \bigl( \lk(X[A],B_1) \bigr) \rightarrow  \ldots ~~~.
\end{equation}
${\bf P(k,m) \Rightarrow P(k+1,m-1):}~~$ Suppose $X$ satisfies
$P(k,m)$ and let $B_1 \subset A_1 \subset V$ such that $|V|-|A_1|
= k+1$ and $|B_1| \leq m-1$. Choose a $v \in V-A_1$ and let $A=A_1
 \cup\{v\}$, $B=B_1\cup \{v\}$. Let $i \geq d$, then by the assumption on $X$,
both the second and the fourth terms in (\ref{mayer}) vanish. It
follows that $\th_i \bigl( \lk(X[A_1],B_1) \bigr)=0$ as required.
\\
${\bf P(k+1,m-1) \Rightarrow P(k,m):}~~$ Suppose $X$ satisfies
$P(k+1,m-1)$ and let $B \subset A \subset V$ such that $|V|-|A|
\leq  k$ and $|B| =m$. Choose a $v \in B$ and let $A_1=A-v$,
$B_1=B-v$. Let $i \geq d$, then by the assumption on $X$, both the
first and the third terms in (\ref{mayer}) vanish. It follows that
$\th_i \bigl( \lk(X[A],B) \bigr)=0$ as required. {\enp} \noindent
{\bf Proof of Proposition \ref{folk}:} Let $X$ be a complex on $n$
vertices. Then (i) is equivalent to $P(n,0)$, while (ii) is
equivalent to $P(0,n)$. On the other hand, $P(n,0)$ and $P(0,n)$
are equivalent by  Claim \ref{st-lk}. {\enp} \noindent {\bf Proof
of Theorem \ref{leray}:} By induction it suffices to consider the
$r=2$ case. Let $X,Y$ be complexes on $V$ with
$\LK(X)=a~,~\LK(Y)=b$, and let $k>a+b$. Then for any $\sigma \in
Y$ and for any $i,j$ such that $i+j=k$, either $i>a$ hence
$\dth_{i-1}(X[\sigma])=0$, or $j>b$  which by Proposition
\ref{folk} implies that $\dth_{j-1}(\lk(Y,\sigma))=0$.  By Theorem
\ref{dimt} it then follows that $\dth_{k-1}(X \cap Y)=0$.
Therefore
\begin{equation}
\label{twoi} \LK(X \cap Y) \leq \max_{S \subset V} \bigl(\LK
(X[S])+\LK (Y[S]) \bigr)=\LK(X)+\LK(Y)~~.
\end{equation}
Next, let $k \geq \LK (X)+\LK (Y)+1$. Then by (\ref{twoi}) and the
Mayer-Vietoris sequence
$$ \rightarrow  \th_k(X) \oplus \th_k(Y)  \rightarrow \th_k(X \cup Y) \rightarrow
\th_{k-1}(X \cap Y) \rightarrow$$ it follows that $\th_k(X \cup Y)
=0$.  Hence
$$\LK (X\cup Y) \leq
\max_{S \subset V}\bigl( \LK (X[S])+\LK (Y[S])+1 \bigr) =
\LK(X)+\LK(Y)+1.$$ {\enp}

\noindent {\bf Acknowledgment:} We thank J\"{u}rgen Herzog for his
comments concerning the commutative algebra aspects of this work.

\end{document}